\documentclass{amsproc}

\usepackage{a4wide}
\usepackage{amsmath}
\usepackage{enumerate}

\usepackage{amsfonts, amssymb,amsthm,amsgen}

\def\gsim{\mathfrak {sim}}
\def\simil{\mathfrak {sim}}
\def\gso{\mathfrak{so}}
\def\so{\mathfrak{so}}
\def\sp{\mathfrak{sp}}
\def\su{\mathfrak{su}}
\def\u{\mathfrak{u}}
\def\spin{\mathfrak{spin}}

\newtheorem{Th}{Theorem}
\newtheorem{Prop}{Proposition}
\newtheorem{Cor}{Corollary}
\newtheorem{Lem}{Lemma}
\newtheorem{Def}{Definition}
\newtheorem{Ex}{Exercise }
\def\bt{\begin{Th}}
\def\et{\end{Th}}
\def\bp{\begin{Prop}}
\def\ep{\end{Prop}}
\def\bc{\begin{Cor}}
\def\ec{\end{Cor}}
\def\bl{\begin{Lem}}
\def\el{\end{Lem}}
\def\bd{\begin{Def}}
\def\ed{\end{Def}}
\def\bex{\begin{Ex}}
\def\eex{\end{Ex}}

\def\be{\begin{equation}}
\def\ee{\end{equation}}
\def\ben{\begin{enumerate}}
 \def\een{\end{enumerate}}
\def\ba{\begin{array}{rlll}}
\def\ea{\end{array}}
\def\bea{\begin{eqnarray}}
\def\eea{\end{eqnarray}}
\def\bean{\begin{eqnarray*}}
\def\eean{\end{eqnarray*}}

\catcode`@=11 \@addtoreset{equation}{section} \catcode`@=12

\def\Real{\mathbb{R}}
\def\Co{\mathbb{C}}

\def\g{\mathfrak{g}}
\def\h{\mathfrak{h}}

\def\p{\partial}
\def\R{\mathcal{R}}
\def\P{\mathcal{P}}

\title{{Classification of third-order symmetric Lorentzian manifolds}}

\author{Anton S. Galaev}

\begin{document}
\maketitle

\begin{abstract} Third-order symmetric Lorentzian manifolds, i.e. Lorentzian manifold with zero third
derivative of the curvature tensor, are classified. These
manifolds are exhausted by a special type of pp-waves, they
generalize Cahen-Wallach spaces and second-order symmetric
Lorentzian spaces.



\end{abstract}

\section{Introduction}

Symmetric pseudo-Riemannian manifolds constitute  an important
class of spaces. A direct generalization of these manifolds is
provided by the so-called $k$-th order symmetric pseudo-Riemannian
spaces $(M,g)$ satisfying the condition $$\nabla^k R=0,\quad
\nabla^{k-1} R\neq 0,$$ where $k\geq 1$ and $R$ is the curvature
tensor of $(M,g)$. For Riemannian manifolds, the condition
$\nabla^k R=0$ implies $\nabla R=0$ \cite{Tanno72}. On the other
hand, there exist pseudo-Riemannian $k$-th order symmetric spaces
with $k\geq 2$, see e.g. \cite{Kai,Sen08}.

The fundamental paper \cite{Sen08} by J.~M.~M.~Senovilla   starts
detailed investigation of symmetric Lorentzian spaces of second
and higher orders. It contains many  interesting results about
such manifolds and their potential physical applications, e.g.
Penrose limit type constructions, usage of super-energy tensors,
higher order Lagrangian theories, supergravity, vanishing of
quantum fluctuations. In particular, it is proven there that any
second-order symmetric Lorentzian space admits a parallel null
vector field, and it is conjectured that this property holds for
symmetric Lorentzian spaces of higher orders.

A classification of four-dimensional second-order symmetric
Lorentzian spaces is obtained by O.~F.~Blanco, M.~S\'anchez,
J.~M.~M.~Senovilla in the paper \cite{Sen10}. The result is based
on the Petrov classification of the conformal Weyl curvature
tensors. In \cite{twosym} D.~V.~Alekseevsky and the author
classified second-order Lorentzian symmetric manifolds. For that,
the methods of the holonomy theory were used. Since the tensor
$\nabla R$ is parallel, its values at points are annihilated by
the holonomy algebra, this algebraic condition allowed to find the
exact form of $\nabla R$; then it was shown that the Weyl
conformal curvature tensor is parallel and the results of
A.~Derdzinski and W.~Roter \cite{DR09} about such spaces were
applied. An alternative proof that used solutions of some PDE
systems was obtained by the authors of \cite{Sen10} in
\cite{Sen13}. Now we know that a second-order Lorentzian symmetric
manifold is locally a product of a locally symmetric Riemannian
manifold and of a Lorentzian manifold with the metric
$$g=2dvdu+\sum_{i=1}^n(dx^i)^2+H(du)^2,\quad H=(H_{1ij} u+H_{0ij})x^ix^j,$$
where $H_{1ij}$  and $H_{0ij}$ are  symmetric real matrices. In
\cite{Sen13} it is shown also that a simply connected geodesically
complete second-order Lorentzian symmetric manifold is a global
product of a (possibly trivial) Riemannian symmetric manifold and
of $\Real^{n+2}$ with the above metric.

This paper is motivated by the lectures of M.~S\'anchez and
J.~M.~Senovilla \cite{San13,Sen14}, where the problems of
classification of the higher order, and first of all of the third
order, symmetric Lorentzian manifolds are discussed. In the
present paper we classify third-order Lorentzian symmetric spaces.
The main result can be stated as follows.

\begin{Th}\label{ThMain} Let $(M,g)$ be a locally indecomposable Lorentzian manifold of dimension
$n+2\geq 4$. Then $(M,g)$ is third-order symmetric if and only if
locally there exist coordinates $v,x^1,...,x^n,u$ such that
$$g=2dvdu+\sum_{i=1}^n(dx^i)^2+H(du)^2,\quad H=(H_{2ij} u^2+ H_{1ij} u+H_{0ij})x^ix^j,$$
where   $H_{2ij}$, $H_{1ij}$ and $H_{0ij}$ are symmetric real
matrices, the matrix $H_{2ij}$ is nonzero, and it can be assumed
to be diagonal.
\end{Th}

By the Wu Theorem \cite{Wu}, any Lorentzian manifold $(M,g)$ is
either locally indecomposable, or it is locally a product of a
Riemannian manifold $(M_1,g_1)$, and of a locally indecomposable
Lorentzian manifold $(M_2,g_2)$. The manifold $(M,g)$ is
third-order symmetric if and only if $(M_1,g_1)$ is locally
symmetric and $(M_2,g_2)$ is third-order symmetric. Consequently,
Theorem \ref{ThMain} provides the complete local classification of
third-order symmetric Lorentzian manifolds.

For the proof of Theorem \ref{ThMain}, we use extend the ideas
from \cite{twosym}. The assumption that a Lorentzian manifold
$(M,g)$ is third-order symmetric implies that the holonomy algebra
of $(M,g)$ at a point $x\in M$ annihilates the tensor $\nabla^2
R_x\neq 0$. This allows to find that exact form of $\nabla^2 R$.
Unlike \cite{twosym} we do not use the Weyl tensor, but find other
tricks that allow to show that the manifold is a pp-wave. Then the
condition $\nabla^3 R=0$ and simple computations allow us to find
the coordinate form of the metric.

In particular, we prove the following two theorems conjectured by
J.~M.~M. Senovilla \cite{Sen10,Sen13,San13} for higher-order
symmetric Lorentzian manifolds.

   \bt\label{Thso} Let $(M,g)$ be Lorentzian manifold of dimension $n+2$ with the holonomy algebra
$\so(1,n+1)$ and such that $\nabla^3R=0$. Then $(M,g)$ is locally
symmetric.\et

   \bt\label{Thparal} Any simply connected third-order symmetric Lorentzian manifold $(M,g)$
    admits a parallel null
vector field.\et

In \cite{Sen13,San13} it is shown that if a complete
simply-connected Lorentzian manifold $(M,g)$ is locally isometric
to the product of a Riemannian symmetric space and a Lorentzian
space with the metric as in Theorem \ref{ThMain}, then $(M, g)$ is
globally isometric to one of such products. This  implies

\begin{Th} Let $(M,g)$ be a simply connected
geodesically complete third-order symmetric Lorentzian
 manifold. Then $(M,g)$ is a product of a (possibly trivial)
 symmetric Riemannian manifold and of $\Real^{n+2}$ with the
 metric from Theorem \ref{ThMain}.
\end{Th}

Finally we discuss possible extension of these results to the case
of higher order symmetric Lorentzian manifolds.

\section{Holonomy algebras of Lorentzian manifolds}
 We recall some basic  facts  about the holonomy groups of  Lorentzian manifolds that can be found in
 \cite{BB-I,ESI,Leistner}.
Let $(M,g)$ be a  Lorentzian $(n+2)$-dimensional  manifold  and
$\g\subset\so(1,n+1)$ be its holonomy algebra at a point $x\in M$,
i.e. the Lie algebra of the holonomy group at that point. Denote
the tangent space $T_xM$ by $V$ and the metric $g_x$ simply by
$g$.

  The  manifold $(M,g)$ is locally indecomposable (i.e. locally is not a direct product
  of two pseudo-Riemannian manifolds) if and only if  the  holonomy
  algebra $\g\subset\so(1,n+1)$ is weakly irreducible, i.e. it  does not preserve  any
   proper nondegenerate  subspace of the tangent space.

 Any  weakly irreducible  holonomy algebra  $\g\subset\so(1,n+1)$ different from  the
    Lorentz Lie algebra $\so(1,n+1)$ preserves a null line
$\Real p$ of the tangent space. Denote by $\simil(n)$ the maximal
subalgebra of $\so(1,n+1)$ preserving $\Real p$. The Lie algebra
$\so(1,n+1)$ will be identified with the space of bivectors
$\Lambda^2 V$ in such a way that $$(X\wedge
Y)Z=g(X,Z)Y-g(Y,Z)X,\quad X,Y,Z\in V.$$ Choose any null vector
$q\in V$ such that $g(p,q)=1$. Let $E\subset V$ be the Euclidean
subspace orthogonal to $p$ and $q$. We get the decomposition
\be\label{decV} V=\Real p\oplus E\oplus \Real q.\ee Let
$e_1,...,e_n$ be an orthonormal basis in $E$. We get the
decomposition into the direct sum of vector subspaces
$$\simil(n)=\Real p\wedge q+\so(n)+p\wedge E,$$
where $\so(n)=\so(E)=\Lambda^2 E$. If $\g\subset\simil(n)$ is an
arbitrary subalgebra, then the  $\gso(n)$-projection of $\g$ is
called {\it the orthogonal part of} $\g$. Decomposition
\eqref{decV} is a $|1|$-grading of $V$ with the grading element
$p\wedge q$.

The weakly irreducible Lorentzian holonomy algebras
$\mathfrak{g}\subset\gsim(n)$ are the following:

\begin{tabular}{llll} (type 1) & $\Real p\wedge
q+\h+p\wedge E$, & (type 2) & $\h+p\wedge E$,\\
(type 3) & $\{\varphi(A)p\wedge q+A|A\in\h\}+p\wedge E$,& (type 4)
& $\{A+p\wedge \psi(A) |A\in\h\}+p\wedge E_1$,\end{tabular}

where $\h\subset\gso(n)$ is a Riemannian holonomy algebra;
$\varphi:\h\to\Real$ is a non-zero linear map that is zero on the
commutant~$\h'=[\h,\h]$; for the last algebra, $E=E_1\oplus E_2$
is an orthogonal decomposition, $\h$ annihilates $E_2$, i.e.
$\h\subset \gso(E_1)$, and $\psi:\h\to E_2$ is a surjective linear
map that is zero on the commutant $\h'$.

A locally indecomposable simply connected Lorentzian manifold
admits a parallel null vector field if and only if its holonomy
group is of type 2 or 4.

Let $\g\subset\gsim(n)$ be the holonomy algebra of the Lorentzian
manifold $(M,g)$ and $\h\subset\gso(E)$ be its orthogonal part.
Then there exist the decompositions
\begin{equation}\label{LM0A}E=E_0\oplus E_1\oplus\cdots\oplus E_r,\quad
\h=\{0\}\oplus\h_1\oplus\cdots\oplus\h_r\end{equation} such that
$\h$ annihilates $E_0$, $\h_i(E_j)=0$ for $i\neq j$, and
$\h_i\subset\gso(E_i)$ is an irreducible subalgebra for $1\leq
i\leq r$.

  \section{The  holonomy algebra of a third-order symmetric Lorentzian manifold}
 A pseudo-Riemannian manifold $(M,g)$ with the curvature tensor $R$  is called a $k$-symmetric space  if
   $$ \nabla^k R = 0,\quad   \nabla^{k-1}R \neq 0.  $$

   So,
    one-symmetric  spaces  are  the same as nonflat locally  symmetric   spaces ($\nabla R =0$, $R\neq 0$).

Remark that for a Riemannian manifold the condition $\nabla^kR=0$
implies $\nabla R=0$ \cite{Tanno72}.

    All indecomposable  simply connected Lorentzian symmetric spaces are exhausted by
     the De Sitter  and the anti De Sitter spaces and
    the  Cahen-Wallach spaces. The last spaces have  the commutative holonomy  algebra $p\wedge E$.

Below we will prove the following

  \bt\label{Thhol}  The  holonomy  algebra of an $(n+2)$-dimensional
  locally indecomposable third-order symmetric Lorentzian manifold $(M,g)$  is
  $p \wedge E \subset \simil(V)$. \et

It is known that any $(n+2)$-dimensional  Lorentzian manifold with
the holonomy algebra $p\wedge E$ is a pp-wave (see e.g. \cite[Sect
5.4]{ESI}), i.e. locally there exist coordinates $v,x^1,...,x^n,u$
such that the metric $g$ can be written in the form
$$g=2dvdu+\sum_{i}(dx^i)^2+H(du)^2,\quad \p_v H=0.$$
We will need only to decide which functions $H$ correspond to the
third-order symmetric spaces.

\section{Algebraic curvature tensors}\label{secR}

For a subalgebra $\g\subset\gso(V)$ define {\it the space of
algebraic curvature tensors of type} $\mathfrak{g}$,
$$\R(\mathfrak{g})=\{R\in \Lambda^2 V^*\otimes\mathfrak{g}\,|\, R(u,v)w+R(v, w)u+R(w, u)v=0
\text{ for all } u,v,w\in V\}.$$  If $\mathfrak{g}\subset\gso(V)$
is the holonomy algebra of   a  manifold $(M,g)$, where $V=T_xM$
is tangent space at some point $x\in M$,
 then the curvature tensor $R_x$ of $(M,g)$ belongs to
$\R(\mathfrak{g})$. The spaces $\R(\g)$ for holonomy algebras of
Lorentzian manifolds are found in \cite{Gal1,onecomp}. For
example, let $\g=\Real p\wedge q+\h+p\wedge E$. For a subalgebra
$\h\subset\gso(n)$ define the space
$$\mathcal{P}(\h)=\{P\in E^*\otimes\h\,|\,g(P(X)Y,Z)+g(P(Y)Z,X)+g(P(Z)X,Y)=0\text{ for all }
X,Y,Z\in E\}.$$ Any $R\in\R(\g)$ (considered as a tensor of type
$(4,0)$) can be written as the sum
$$R=R_0+P+T+v+L,$$
where \begin{align*} R_0&=R_0^{ijkl}(e_i\wedge
e_j)\odot(e_k\wedge e_l)\in\R(\h),\\
P&=P^{ijk}(e_i\wedge
e_j)\odot(p\wedge e_k),\quad P^{ijk}\in\Real,\quad P(\cdot,q)\in\P(\h),\quad P^{ijk}+P^{jki}+P^{kij}=0,\\
T&=T^{ij}(p\wedge e_i)\odot(p\wedge e_j),\quad T^{ij}\in\Real,\quad T^{ij}=T^{ji},\\
v&=v^i(p\wedge q)\odot(p\wedge e_i), \quad v^i\in\Real,\\
L&=\lambda (p\wedge q)\odot(p\wedge q),\quad \lambda \in
\Real.\end{align*} Here for bivectors $\omega$ and $\theta$, we
write
$$\omega\odot\theta=\omega\otimes\theta+\theta\otimes\omega.$$
The decomposition \eqref{LM0A} implies
$$R_0=R_{01}+\cdots+R_{0r},\quad P=P_1+\cdots+P_r,\quad
R_{0\alpha}\in\R(\h_\alpha),\quad
P_{\alpha}(\cdot,q)\in\P(\h_\alpha).$$

{\bf Notation.} If $S\in \otimes^r V\,\otimes\R(\g)$, then we
write
$$S=e_{a_1}\otimes\cdots\otimes e_{a_r}\otimes R^{a_1\cdots
a_r},\quad R^{a_1\cdots a_r}\in\R(\g),$$ where we assume that the
indices take the values $p,1,...,n,q$, and that $e_p=p$, and
$e_q=q$. Next, we write $$R^{a_1\cdots a_r}=R_0^{a_1\cdots
a_r}+P^{a_1\cdots a_r}+T^{a_1\cdots a_r}+v^{a_1\cdots
a_r}+L^{a_1\cdots a_r},$$ where e.g. $$P^{a_1\cdots
a_r}=P^{a_1\cdots a_rijk}(e_i\wedge e_j)\odot(p\wedge e_k),\quad
T^{a_1\cdots a_r}=T^{a_1\cdots a_rij}(p\wedge e_i)\odot(p\wedge
e_j).$$

\vskip0.2cm

 Now we define the space of algebraic covariant derivatives of the curvature
 tensors
$$\nabla\R(\mathfrak{g})=\{S\in V^*\otimes\R(\mathfrak{g})\,|\, S_u(v,w)+S_v(w,u)+S_w(u,v)=0
\text{ for all } u,v,w\in V\}.$$ If $\mathfrak{g}\subset\gso(V)$
is the holonomy algebra of a  manifold $(M,g)$ at a point $x\in
M$, then $\nabla R_x\in\nabla\R(\mathfrak{g})$. The decomposition
of the space $\nabla\R(\gso(r,s))$ into irreducible
$\gso(r,s)$-modules is found in \cite{Str88}, see also
\cite{nablaR}.

Let us find the space $\nabla\R(\g)$ for $\g=\Real p\wedge
q+\h+p\wedge E$.

\bt\label{ThnabR} Any $S\in\nabla\R(\g)$ has the form
$$S=p\otimes R^p+e_t\otimes R^t+q\otimes R^q,\quad R^p,R^t,R^q\in \R(\g),$$
and it holds $$P^{pijk}=T^{ijk}-T^{jik},\quad
P^{tijk}=2R_0^{ijtk},\quad R_0^q=0,\quad P^q=0,$$
$$e_t\otimes R^t_0\in\nabla\R(\h),\quad v^{ij}=2T^{qij},\quad v^{qi}=2\lambda^i.$$
\et

{\bf Proof.} We may write $S=p\otimes R^p+e_t\otimes R^t+q\otimes
R^q$ for some elements $R^p,R^t,R^q\in \R(\g)$. The equality
$$S_p(e_t,e_s)+S_{e_t}(e_s,p)+S_{e_s}(p,e_t)=0$$ can be rewritten
in the form $$4R_0^{qijts} e_i\wedge e_j+2P^{qtsk}p\wedge e_k=0.$$
This implies $R_0^q=0$ and $P^q=0$. Considering the vectors
$e_m,e_s,e_t$, we get $e_t\otimes R^t_0\in\nabla\R(\h)$, and
$$P^{msrk}+P^{srmk}+P^{rmsk}=0.$$ this means that
$P^{msrk}e_m\otimes e_k\otimes (e_s\wedge e_r)\in\R(\h)$, and, in
particular, $P^{msrk}=-P^{ksrm}$. Considering the vectors $p$,
$e_m$, $q$, we get $$-2T^{qim}p\wedge e_i-v^{qm}p\wedge
q+2\lambda^m p\wedge q+v^{mi}p\wedge e_i=0.$$ Consequently,
$v^{mi}=2T^{qim}$, and $v^{qm}=2\lambda^m$. Using the vectors $q$,
$e_t$, $e_s$, we get the rest of the equalities.
  \qed

\section{Proof of Theorem \ref{Thso}}

Let $(M,g)$ be an $(n+2)$-dimensional Lorentzian manifold with the
holonomy algebra $\g=\so(1,n+1)$ and the property $\nabla^3R=0$.
It is noted in \cite{Sen10,San13} that using the methods from
\cite{Tanno72} it can be shown that $$g(\nabla^2R,\nabla^2R)=0.$$
The tensor $\nabla^2 R$ is
 parallel and annihilated by the holonomy algebra.
Consequently, for each point $x\in M$,
$$\nabla^2R_x:T_xM=V\to\nabla\R(\g),\quad X\mapsto (\nabla_X\nabla R)_x$$
is an $\g$-equivariant map. The multiplicity of the  $\g$-module
$V$ in the space $\nabla\R(\g)$ is one \cite{Str88,nablaR}.
Consequently, the multiplicity of the trivial $\g$-module $\Real$
in $V\otimes\nabla\R(\g)$ is one as well, and $\nabla^2R_x$
belongs to this submodule. The extension of the metric $g_x$ to
the space $V\otimes\nabla\R(\g)$ is non-degenerate and the space
$V\otimes\nabla\R(\g)$ can be decomposed into an orthogonal direct
sum of $\g$-invariant modules. Hence the restriction of $g_x$ to
$\Real\subset V\otimes\nabla\R(\g)$ is non-degenerate. We conclude
that $\nabla^2 R_x=0$, i.e. $\nabla^2R=0$. Results of
\cite{twosym,Sen08} imply that $\nabla R=0$. \qed

\section{Walker coordinates and a reduction lemma}\label{seccoord}

Let $(M,g)$ be a locally indecomposable third-order Lorentzian
manifold with the weakly irreducible holonomy algebra
$\g\subset\so(1,n+1)$. From Theorem \ref{Thso} it follows that
$\g\neq\so(1,n+1)$. Since $\gso(1,n+1)$ is the only irreducible
holonomy algebra \cite{ESI}, it follows that $\g\subset\gsim(n)$.

Let $(M,g)$ be a Lorentzian manifold  with the holonomy algebra
$\g\subset\gsim(n)$. Then $(M,g)$ admits (locally) a parallel
distribution of null lines.  According to \cite{Walker}, locally
there exist the so called Walker coordinates $v,x^1,...,x^n,u$
such that the metric $g$ has the form
\begin{equation}\label{Walker} g=2dvdu+h+2Adu+H (d
u)^2,\end{equation} where $h=h_{ij}(x^1,...,x^n,u)d x^id x^j$ is
an $u$-dependent family of Riemannian metrics,  $A=A_i(x^1,
\ldots, x^n,u)d x^i$ is an $u$-dependent family of one-forms, and
$H=H(v,x^1,...,x^n,u)$ is a local function on $M$. Consider the
local frame \be \label{frame}p=\p_v,\quad X_i=\p_i-A_i\p_v,\quad
q=\p_u-\frac{1}{2}H\p_v.\ee Let $E$ be the distribution generated
by the vector fields $X_1$,...,$X_n$. Clearly, the vector fields
$p$, $q$ are null, $g(p,q)=1$, the restriction of $g$ to $E$ is
positive definite, and $E$ is orthogonal to $p$ and $q$. The
vector field $p$ defines the parallel distribution of null lines
and it is recurrent, i.e. $\nabla p=\theta\otimes p$, where
$\theta=\frac{1}{2}\p_vHdu$. Since the manifold is locally
indecomposable, any other recurrent vector field is proportional
to $p$. Next, $p$ is proportional to a parallel vector field if
and only if $d\theta =0$, which is equivalent to
$\p^2_vH=\p_i\p_vH=0$. In the last case the coordinates can be
chosen in such a way that $\p_v H=0$ and $\nabla p=\nabla \p_v=0$,
see e.g. \cite{ESI}.

Consider the metric \eqref{Walker}, the vector fields
\eqref{frame} and an orthogonal frame $e_1,...,e_n$ of the
distribution $E$. Then the curvature tensor $R$ of the metric and
its covariant derivatives can be written as in Section \ref{secR}
above with respect to the frame $p,e_1,...,e_n,q$ and all
coefficients being functions.

In \cite{GL10} using results from\cite{Boubel} it is shown that
 there exist Walker coordinates
$$v,\,x_0=(x_0^1,\ldots ,x_0^{n_0}),\ldots ,x_r=(x_r^1,...,x_r^{n_r}),\,u$$
adapted to the decomposition \eqref{LM0A} and in addition  with
the property $A=0$. This means that $$h=h_0+h_1+\cdots+h_r,\quad
h_0=\sum_{i=1}^{n_0}(d x_0^i)^2,\quad
h_\alpha=\sum_{i,j=1}^{n_\alpha}h_{\alpha ij} d x_{\alpha}^id
x_{\alpha}^j,\quad \frac{\p}{\p {x^k_\beta}}h_{\alpha ij}=0 \text{
if }\beta\neq\alpha.$$

For $\alpha=0,...,r$, consider  the submanifolds
$M_{\alpha}\subset M$ defined by $x_\beta=c_\beta$,
$\alpha\neq\beta$, where $c_\beta$ are constant vectors. Then the
induced  metric is given by
$$g_\alpha=2dvdu+h_\alpha+H_\alpha(du)^2.$$

 The proof of the following lemma is the same as one of Lemma 1 from \cite{twosym}.

\begin{Lem}\label{Lemgi} The submanifold $M_\alpha\subset M$ is
totally geodesic. The orthogonal part
 of the holonomy algebra
$\g_{\alpha}$ of the metrics $g_\alpha$ coincides with
$\h_\alpha\subset\gso(E_\alpha)$, which is irreducible for
$\alpha=1,...,r$. If the metric $g$ is third-order symmetric, then
the curvature tensor of each metric $g_\alpha$ satisfies
$\nabla^3R=0$.
 \end{Lem}

Remark that the metric $g_\alpha$ must not be indecomposable.

\section{Proof of Theorem \ref{Thparal}}

We will show that  there are no Lorentzian manifolds with the
property $\nabla^3R=0$ and with the holonomy algebras of type 1 or
3.

\bl\label{Lt1} Let $\g=\Real p\wedge q+\h+p\wedge E$ with no
assumption on $\h\subset\so(n)$, then the subspace of
$V\otimes\nabla\R(\g)$ annihilated by $\g$ is trivial.\el

{\bf Proof.} Let $S\in V\otimes\nabla\R(\g)$ and suppose that it
is annihilated by $\g$. Let us write $$S=p\otimes S^p+e_t\otimes
S^t+q\otimes S^q,$$  where each element $S^p,S^t,S^q\in
\nabla\R(\g)$ is as in Theorem \ref{ThnabR}, e.g.
$$S^p=p\otimes R^{pp}+e_t\otimes
R^{pt}+q\otimes R^{pq},\quad R^{pp},R^{pt},R^{pq}\in \R(\g).$$
Note that $$(p\wedge q)p=-p,\quad (p\wedge q)q=q, \quad (p\wedge
q)e_i=0.$$ Consequently, if $$Q=Q^{a_1\cdots
a_r}e_{a_1}\otimes\cdots\otimes e_{a_r},\quad Q^{a_1\cdots
a_r}\in\Real$$ is a tensor annihilated by $p\wedge q$, then
$Q^{a_1\cdots a_r}=0$ whenever in $Q^{a_1\cdots a_r}$ the number
of indices equals to p is different from the number of indices
equal to $q$. This implies that
$$S^p=q\otimes\lambda^{pq}(p\wedge q)\odot (p\wedge q).$$ The
condition that $p\wedge e_s\in\g$ annihilates $S$ implies the
equations \begin{equation}\label{SSS}S^s=(p\wedge e_s)\cdot
S^p,\quad (p\wedge e_s)\cdot S^t=-\delta_{st}S^q,\quad (p\wedge
e_s)\cdot S^q=0.\end{equation} Using this we obtain
$$S^s=e_s\otimes\lambda^{pq}(p\wedge q)\odot (p\wedge
q)+2q\wedge\lambda^{pq}(p\wedge e_s)\odot (p\wedge q),$$
\begin{multline*}-\delta_{st}S^q=(p\wedge
e_t)\cdot S^s=-\delta_{st}p\otimes\lambda^{pq}(p\wedge q)\odot
(p\wedge q)+2e_s\otimes\lambda^{pq}(p\wedge e_t)\odot (p\wedge
q)\\+2e_t\otimes\lambda^{pq}(p\wedge e_s)\odot (p\wedge
q)+2q\otimes\lambda^{pq}(p\wedge e_s)\odot (p\wedge
e_t).\end{multline*} Taking $s\neq t$, we get $\lambda^{pq}=0$.
Consequently, $S=0$. \qed

\begin{Lem}\label{Jinv} Let $J$ be a complex structure on $\Real^{2m}$. Then
the eigenvalues of $J$ on $\odot^2\Real^{2m}$ are zero and pure
imaginary, and the eigenvalues of $J$ on $\odot^3\Real^{2m}$ are
 pure imaginary. \end{Lem}

 {\bf Proof.} As usual, the complexification $\Real^{2m}\otimes\Co=\Co^{2m}$ can be
decomposed as $\Co^{2m}=W\oplus\bar W$, where $W$ is the
eigenspace of the extension of $J$ with eigenvalue $i$, and $\bar
W$ is the eigenspace of the extension of $J$ with eigenvalue $-i$.
Next,
\begin{equation}\label{J2}(\odot^2\Real^{2m})\otimes\Co=\odot^2\Co^{2m}=(\odot^2
W)\oplus(W\otimes \bar W)\oplus(\odot^2\bar W)).\end{equation}
This shows that the eigenvalues of $J$ on $\odot^2\Real^{2m}$ are
$2i$, $0$ and $-2i$. Similarly,
$$(\odot^3\Real^{2m})\otimes\Co=\odot^3\Co^{2m}=(\odot^3
W)\oplus(\odot^2 W\otimes \bar W)\oplus ( W\otimes \odot^2\bar
W)\oplus(\odot^3\bar W)),$$ i.e. the eigenvalues of $J$ on
$\odot^2\Real^{2m}$ are $3i$, $i$, $-i$ and $-3i$. \qed

\bl\label{Lt3} Let $\g\subset\simil(n)$ be of type 3 with an
irreducible orthogonal part $\h\subset\so(n)$, then the subspace
of $V\otimes\nabla\R(\g)$ annihilated by $\g$ is trivial.\el

{\bf Proof.} Since $\g$ is of type 3, and $\h\subset\so(n)$ is
irreducible, it holds $\h\subset\mathfrak{u}(m)$, $n=2m$, and
$$\g=\Real(p\wedge q+cJ)+\h'+p\wedge E,\quad c\in\Real,\quad c\neq 0,$$
where $J$ is the complex structure on $\Real^{2m}$. Suppose that
$S\in V\otimes\nabla\R(\g)$ is annihilated by $\g$. Let as in
Lemma \ref{Lt1} $$S=p\otimes S^p+e_t\otimes S^t+q\otimes S^q,\quad
S^p,S^t,S^q\in \nabla\R(\g).$$ Since $\g$ is of type 3, it holds
$\lambda^{ab}=0$. Let $\xi=p\wedge q+cJ\in\g$. It is clear that
$\xi\cdot (p\otimes S^p)=0$. Consequently,
$$S^p=\xi\cdot S^p=-p\otimes R^{pp}+p\otimes\xi\cdot R^{pp}+\xi\cdot (e_t\otimes
R^{pt})+q\otimes R^{pq}+q\otimes\xi\cdot R^{pq}.$$ We get the
equations
$$\xi\cdot R^{pp}=2R^{pp},\quad \xi\cdot(e_t\otimes R^{pt})=e_t\otimes
R^{pt},\quad \xi\cdot R^{pq}=0.$$ Note that for each $R\in\R(\g)$
it holds $$R=R_0+P+T+v,\quad \xi\cdot R= -P+cJP-2T+cJT-v+cJv.$$
Using this we get $$R_0^{pp}=0,\quad cJP^{pp}=3P^{pp},\quad
cJT^{pp}=4T^{pp},\quad cJv^{pp}=3v^{pp},$$
$$cJ(e_t\otimes R_0^{pt})=e_t\otimes R_0^{pt},\quad CJ(e_t\otimes
T^{pt})=3e_t\otimes T^{pt}.$$ Note that since
$\h\subset\mathfrak{u}(m)$, it holds $JR_0=0$ for each
$R_0\in\R(\h)$, and $JP=P^{ijk}(e_i\wedge e_j)\otimes J e_k$ for
each $P=P^{ijk}(e_i\wedge e_j)\otimes e_k\in \P(\h)$. We conclude
that
$$R_0^{pt}=R_0^{pp}=P^{pp}=T^{pp}=v^{pp}=0,$$
in particular, $R^{pp}=0$ (the equality $T^{pp}=0$ follows from
the previous lemma). From Theorem \ref{ThnabR} it follows that
$P^{pt}=0$, and the tensor $T^{pijk}$ is symmetric in $i,j,k$.
From the equality $cJ(e_t\otimes T^{pt})=3e_t\otimes T^{pt}$ and
the previous lemma it follows that $T^{pt}=0$. We get that
$R^{pt}=v^{pt}$. The equality $\xi\cdot S^{pq}=0$ implies
$S^{pq}=0$. From this and Theorem \ref{ThnabR} it follows that
$v^{pt}=0$. Thus, $S^p=0$. From \eqref{SSS} it follows that $S=0$.
\qed

Let $(M,g)$ be simply connected third-order symmetric Lorentzian
manifold. Suppose that it does not admit a parallel null vector
field. Then its holonomy algebra is either of type 1 or 3. Recall
that there are no neither second-order symmetric nor locally
symmetric Lorentzian manifolds with holonomy algebras of type 1 or
3. If the holonomy algebra is of type 1, then from Lemma \ref{Lt1}
it follows that $\nabla^2R=0$, which gives a contradiction. Hence
the holonomy algebra is of type 3. Similarly, from Lemma \ref{Lt3}
it follows that the orthogonal part $\h\subset\so(n)$ is not
irreducible. Consider the metrics $g_\alpha$ as in Section
\ref{seccoord}. Since $\g$ is of type 3, it holds $\p^2_vH=0$ and
$\p_v\p_{x^{i_\alpha}_\alpha}H_\alpha\neq 0$ for some $\alpha$ and
$i_\alpha$ \cite{comphol}. This implies that the metric $g_\alpha$
is indecomposable with $\nabla^3 R=0$ and its holonomy algebra is
of type 3 with irreducible orthogonal part
$\h_\alpha\subset\so(n_\alpha)$. From Lemma \ref{Lt3} it follows
that $\nabla^2R=0$, this gives a contradiction. \qed

\section{Proof of Theorem \ref{Thhol}}

Now we assume that $(M,g)$ is an indecomposable third-order
symmetric Lorentzian manifold with the holonomy algebra
$\g\subset\simil(n)$ annihilating the null vector $p\in V$.  The
metric $g$ written as \eqref{Walker} than satisfies $\p_v H=0$.
Clearly it holds $\Gamma_{va}^b=0$, $R_{vabc}=0$, and $\nabla
R_{abcd;e}=0$ whenever on of the indices is $v$. In particular,
for any $x\in M$ it holds $\nabla^2R_x\in (\Real p\oplus
E)\otimes\nabla\R(\g)$.

\bl\label{Lt2} Let $\g=\h+p\wedge E$ with $\h\subset\so(n)$ being
an arbitrary subalgebra, then any element $S\in (\Real p\oplus
E)\otimes\nabla\R(\g)$ annihilated by $\g$ is of the form
$$S=(T^{pij} p\otimes p+ T^{pkij}p\otimes e_k-T^{pkij}e_k\otimes p)
\otimes (p\wedge e_i)\odot (p\wedge e_j),$$ such that  the tensors
$T^{pij}e_i\otimes e_j$ and $T^{pkij}e_k\otimes e_i\otimes e_j$
are symmetric and annihilated by $\h$.\el

{\bf Proof.} Let $S\in (\Real p\oplus E)\otimes\nabla\R(\g)$ and
suppose that it is annihilated by $\g$. Let us write $$S=p\otimes
S^p+e_t\otimes S^t,$$  where the elements $S^p,S^t\in
\nabla\R(\g)$ are as in Theorem \ref{ThnabR}. Since $\g$ is of
type 2, all tensors $L$ and $v$ are zero. Let $A\in \h$. Then it
is clear that $A\cdot S^p=0$. Since
$$S^p=p\otimes (R_0^{pp}+P^{pp}+T^{pp})+e_m\otimes
(R_0^{pm}+P^{pm}+T^{pm}),$$ we conclude that
$$A\cdot P^{pp}=0,\quad A\cdot T^{pp}=0,\quad A\cdot (e_m\otimes
R_0^{pm})=0,\quad A\cdot (e_m\otimes T^{pm})=0.$$ The spaces
$\P(\h)$ do not contain elements annihilated by $\h$
\cite{onecomp}, consequently, $P^{pp}=0$. From this and Theorem
\ref{ThnabR} it follows that the tensor $T^{pmij}$ is symmetric.
The spaces $\R(\h)$ do not contain submodules isomorphic to $E$
\cite{Al}, this implies that $e_m\otimes R_0^{pm}=0$, i.e.
$R_0^{pm}=0$.

The condition that $p\wedge e_s\in\g$ annihilates $S$ implies the
equations \begin{equation}\label{SSS1}S^s=(p\wedge e_s)\cdot
S^p,\quad (p\wedge e_s)\cdot S^t=0.\end{equation} Using this we
obtain
$$S^s=4p\otimes R_0^{ppijks}(e_i\wedge e_j)\odot(p\wedge
e_k)-p\otimes P^{ps}-p\otimes T^{ps}+2e_t\otimes P^{ptisk}(p\wedge
e_i)\odot(p\wedge e_k).$$ Next, $$0=(p\wedge e_m)\cdot
S^s=p\otimes(8 R_0^{ppimks}-2P^{psimk}-2P^{pmisk})(p\wedge
e_i)\odot(p\wedge e_k).$$ From Theorem \ref{ThnabR} it follows
that $P^{pmisk}=2R_0^{ppismk}$. This and the last equality imply
that $R_0^{pp}=0$ and $P^{pt}=0$. This proves the lemma. \qed

Now we continue the proof of Theorem \ref{Thhol}.

Consider the metric \eqref{Walker}, the vector fields
\eqref{frame} and an orthogonal frame $e_1,...,e_n$ of the
distribution $E$. Then $\nabla^2 R$ can be written as in Lemma
\ref{Lt2} with respect to the frame $p,e_1,...,e_n,q$ and the
elements $T^{pij}$ and $T^{pkij}$ being functions. Note that the
1-form $du$ is dual to the vector field $p$ and it is parallel.

Suppose that $T^{pkij}=0$. This means that
$$\nabla^2 R=du\otimes du\otimes T^{pp},$$ where
$\nabla^2 R$ is considered as the tensor of type $(4,2)$. Since
$\nabla^3R=0$ and $\nabla du=0$, we get  $\nabla T^{pp}=0$.
Consequently, $$\nabla^2 R=\nabla^2\frac{u^2}{2} T^{pp},$$ i.e.
$$\nabla\left(\nabla R-\nabla \frac{u^2}{2} T^{pp}\right)=0.$$
We see that value of the tensor $\nabla R-\nabla
\frac{u^2}{2}T^{pp}$ at any point $x\in M$  is annihilated by the
holonomy algebra and belongs to the space $\nabla\R(\g)$. As in
\cite[Lemma 3]{twosym} or as in Lemma \ref{Lt2} it can be shown
that
$$\nabla R-\nabla \frac{u^2}{2} T^{pp}=du\otimes T_1,\quad
T_1=T_1^{ij}(p\wedge e_i)\odot (p\wedge e_j).$$ Consequently,
$$\nabla\left(R-\frac{u^2}{2} T^{pp}-uT_1\right)=0.$$ This easily implies
that $$R=\frac{u^2}{2} T^{pp}+uT_1+T_0,\quad T_0=T_0^{ij}(p\wedge
e_i)\odot (p\wedge e_j),$$ i.e. the curvature tensor $R$ is the
same of a pp-wave. We conclude that the metric $g$ is locally a
pp-wave metric, and $\g=p\wedge E$.

Suppose now that $T^{pkij}\neq 0$ and suppose that $\h\neq 0$.
 Consider two cases.

{\bf Case 1.} $\h\subset\so(n)$ is irreducible. The Riemannian
holonomy algebra $\h\subset\so(n)$ annihilates a symmetric
3-tensor. It is known \cite{Nur08} that this happens only for the
irreducible representations of the Lie algebras $\so(3)$,
$\su(3)$, $\sp(3)$, and $F_4$ in dimensions $5$, $8$, $14$, and
$26$, respectively. We do not need such a strong statement and we
prove a weaker one in order to make the exposition more
self-contained. Recall that a symmetric Berger algebra
$\h\subset\so(n)$ is the holonomy algebra of  a symmetric
Riemannian manifold different from $\so(n)$, $\u(\frac{n}{2})$ and
$\sp(\frac{n}{4})\oplus\sp(1)$.

\bl If an irreducible Riemannian holonomy algebra
$\h\subset\so(n)$ admits a nonzero invariant symmetric 3-tensor,
then $\h\subset\so(n)$ is a symmetric Berger algebra. \el

{\bf Proof.} We must show, that the Riemannian holonomy algebras
$\so(n)$, $\u(\frac{n}{2})$, $\sp(\frac{n}{4})\oplus\sp(1)$,
$\su(\frac{n}{2})$, $\sp(\frac{n}{4})$, $G_2\subset\so(7)$ and
$\spin(7)\subset\so(8)$ do not admit nonzero invariant symmetric
3-tensors. We claim that in each case the module $\odot^2\Real^n$
does not contain any submodule isomorphic to $\Real^n$. To check
this it is enough to pass to the complexifications and use tables
from \cite{V-O}, where one may find decompositions of the modules
$\odot^2\Co^n$ for all considered representations. This proves the
claim for the Lie algebras $\so(n)$, $G_2$ and $\spin(7)$. For the
proof of the claim for the Lie algebras $\su(\frac{n}{2})$ and
$\sp(\frac{n}{4})$ also the decomposition \eqref{J2} should be
used (and this will imply the proof of the claim for the Lie
algebras $\u(\frac{n}{2})$ and $\sp(\frac{n}{4})\oplus\sp(1)$).
\qed

Since  $T^{pkij}\neq 0$, we conclude that $\h\subset\so(n)$ is a
symmetric Berger algebra.

From Lemma \ref{Lt2} it follows that the curvature tensor of
$(M,g)$ satisfies \be\label{RdotR}R(e_k,q)\cdot
R=\nabla^2_{e_k;q}R-\nabla^2_{q;e_k}R=2T^{pkij}(p\wedge
e_i)\odot(p\wedge e_j).\ee On the other hand, we write $R=R_0+P+T$
 as in Section
\ref{secR} (for $\g$ under consideration it holds $L=v=0$), and we
get
$$R(e_k,q)\cdot R= P(e_k,q)\cdot (R_0+P+T)+T(e_k,q)\cdot
(R_0+P+T).$$ It can be directly checked that
$$T(e_k,q)\cdot R_0=8T^{ik}\delta_{it}R_0^{jlrt}(e_j\wedge e_l)\odot(p\wedge e_r).$$
Since the action of $P(e_k,q)$ preserves the grading of $V$, we
get
$$P(e_k,q)\cdot P+T(e_k,q)\cdot R_0=0,$$
which can be rewritten in the form
$$(P(e_k,q)\cdot P)(Y,q)+8R_0(T^{ik}e_i,Y),\quad Y\in E.$$
For each symmetric Berger algebra $\h\subset\so(n)$, the space
$\R(\h)$ is one-dimensional and is spanned by a tensor $\R_0$;
next,  $\P(\h)\simeq\Real^n$, and each $P\in\P(\h)$ is of the form
$P=R_0(\cdot,X)$ for some $X\in\Real^n$ \cite{onecomp}.

Consequently, in our situation, if $R_0$ is nonzero on an open set
$U$, then there exists a section $X$ of $E$ over $U$ such that
$P(Y,q)=R_0(Y,X)$ for all sections $Y$ of $E$ over $U$. In
\cite{GL10} it is shown that if we consider the new vector field
$$q'=-\frac{1}{2}g(X,X)p+X+q,$$ and the corresponding distribution
$E'$ with the sections $$Y'=-g(Y,X)p+Y,\quad Y\in E,$$ then in new
notations it holds $P=0$. Next, $R_0(T^{ik}e_i,\cdot)=0$. From the
Bianchi identity it follows that $R_0(Y,Z)T^{ik}e_i=0$ for all
$Y,Z\in E$. Since the values of  $R_0(Y,Z)$ at the point $x$
generate $\h$, and $\h\subset\so(n)$ is irreducible, we get that
$T^{ik}=0$. This implies that $R(e_k,q)\cdot R=0$, i.e.
$T^{pijk}=0$, and we get a contradiction.

Suppose that $R_0=0$. Then $P(e_k,q)\cdot P=0$. This implies $P=0$
\cite{onecomp}. We get $R=T$ is the curvature tensor of a pp-wave,
hence $(M,g)$ is a pp-wave and this is a contradiction.

{\bf Case 2.} $\h\subset\so(n)$ is not irreducible. Then we have
the decomposition \eqref{LM0A}. From Lemma \ref{Lemgi} it follows
that each metric $g_\alpha$, $\alpha=1,...,r$, satisfies
$\nabla^3R=0$. If the metric $g_\alpha$ is indecomposable, then
its holonomy algebra is $\h_\alpha+p\wedge E_\alpha$. According to
the previous case, this is impossible. Consequently the metric is
decomposable, i.e. it is the sum of a Riemannian metric and of a
Lorentzian metric. From Lemma \ref{Lemgi} it follows that the
holonomy algebra of the Riemannian part is isomorphic to
$\h_\alpha$. Hence, the Lorentzian part is of dimension 2 and its
holonomy algebra is either trivial or it is isomorphic to
$\so(1,1)$; since there exists a null parallel  vector field, the
Lorentzian part is flat. We conclude that there exists a null
parallel vector field not proportional to $p$. Clearly, it should
be of the form
$$q_\alpha=-\frac{1}{2}g(X_\alpha,X_\alpha)p+X_\alpha+q,\quad
X_\alpha\in E_\alpha.$$ It holds $R_\alpha(q_\alpha,\cdot)=0$.
Hence, $P_\alpha(Y,q)=-R_{\alpha0}(Y,X_\alpha)$ for all $Y\in
E_\alpha$. Considering as above the new vector field
$$q'=-\frac{1}{2}g(X,X)p+X+q,\quad X=-(X_1+\cdots +X_r),$$
we will get $P=0$. The rest of the proof as in the previous case.
The theorem is proved. \qed

\section{Proof of Theorem \ref{ThMain}}

We know now that the holonomy algebra of the manifold $(M,g)$ is
$p\wedge E$, i.e. the manifold is a pp-wave and locally it is
given by the metric
 \be\label{pp-wave}   g = 2 du dv + \sum_{i=1}^n(dx^i)^2  + Hdu^2,
 \ee
where $H$ is a function  of $x^i$ and $u$, see e.g. \cite[Sect.
5.4]{ESI}. We need only to decide for which functions $H$ the
metric is third-order symmetric. All tensors will be considered as
contravariant. It can be shown that with respect to the frame
$p=\p_v,e_i=\p_{x^i},q=\p_u-\frac{1}{2}\p_v$ it holds
$$\nabla p=0,\quad\nabla e_i=\frac{1}{2}H_{,i} p\otimes p,\quad \nabla q=-\sum_i\frac{1}{2}H_{,i} p\otimes
e_i,$$ where the comma denotes the partial derivative. This
implies \cite{twosym} \begin{align*}  R=&
\frac{1}{2}\sum_{i,j}H_{,ij}(p \wedge e_i) \odot (p \wedge
e_j),\quad \nabla
R=\frac{1}{2}\sum_{i,j}\left(\sum_{k}H_{,ijk}e_k\otimes+H_{,iju}p\right)\otimes(p
\wedge e_i)\otimes (p \wedge e_j ),\\ \nabla^2
R=&\sum_{i,j}\left(\left(\frac{1}{2}H_{,ijuu}-\frac{1}{4}\sum_kH_{,k}H_{,ijk}\right)p\otimes
p\right.\\ &+\left.
 \frac{1}{2}\sum_{k}H_{,ijku}( p \odot e_k) +
 \frac{1}{2}\sum_{k,l}H_{,ijkl}(  e_k \otimes e_l)\right) \otimes (p \wedge e_i)\odot (p \wedge
 e_j).\end{align*}
From Lemma \ref{Lt2} it follows that $H_{,ijkl}=H_{,ijku}=0$.
Consequently,
$$H=H_{ijk}x^ix^jx^k+F_{ij}(u)x^ix^j+G_i(u)x^i+K(u),\quad
H_{ijk}\in\Real.$$ Next,
$$\nabla^2R=\sum_{i,j}\left(\frac{1}{2}H_{,ijuu}-\frac{1}{4}\sum_kH_{,k}H_{,ijk}\right)p\otimes
p\otimes (p \wedge e_i)\odot (p \wedge
 e_j),$$ and the equality $\nabla^3R=0$ is equivalent to to the
 equations
 $$\p_{x^l}(2H_{,ijuu}-\sum_kH_{,k}H_{,ijk})=0,\quad \p_u(2H_{,ijuu}-\sum_kH_{,k}H_{,ijk})=0.$$
The first equation implies $\sum_kH_{,kl}H_{,ijk}=0$. Applying
$\p_{x^s}$, we get $\sum_kH_{,kls}H_{,ijk}=0$. For fixed $l$ and
$i$ this equality means that the squire of the symmetric matrix
$(H_{,ijk})_{j,k=1}^n$ is zero, consequently, $H_{,ijk}=0$. We are
left with the equation $H_{,ijuuu}=0$. This implies that
$$F_{ij}(u)=H_{2ij} u^2+ H_{1ij} u+H_{0ij}.$$  The obtained metric has the same curvature tensor as the
metric from the formulation of the theorem, this implies that the
both metrics are isometric. \qed

{\bf Conclusion.} We have extended and improved the methods from
\cite{twosym} in order to classify the third-order symmetric
Lorentzian manifolds. It is natural to ask about the
classification of the $k$-th order  symmetric Lorentzian manifolds
for arbitrary $k\geq 4$ \cite{Sen10,San13,Sen14}. The description
of the holonomy-invariant tensors in $\otimes^{k-1}V\otimes\R(\g)$
becomes much more complicated. This means that the algebraic
approach developed by us must be intensively combined with
geometric and analytic approaches from \cite{Tanno72,Sen13}. It is
not enough to use  the fact that the tensor $\nabla^{k-1}R$ is
holonomy-invariant, but also the fact that this is the $(k-1)$-th
covariant derivative of the curvature tensor must be applied.
Formulas similar to \eqref{RdotR} should be intensively used.

{\bf Acknowledgements.} I am  thankful to Miguel S\'anchez for
taking my attention to the problem of classification of the
third-order symmetric spaces. I am grateful to Dmitry~V.
Alekseevsky for useful discussions.

\end{document}